\documentclass[12pt, reqno]{amsart}
\usepackage{ amsmath,amsthm, amscd, amsfonts, amssymb, graphicx, color}
\usepackage[bookmarksnumbered, colorlinks, plainpages]{hyperref}
\newtheorem{thm}{Theorem}[section]
\newtheorem{cor}[thm]{Corollary}
\newtheorem{lem}[thm]{Lemma}

\newtheorem{exam}[thm]{Example}
\numberwithin{equation}{section}

\begin{document}

\title{New representation of Drazin inverse for a block complex matrix}

\author{Huanyin Chen}
\author{Marjan Sheibani$^*$}
\address{
Department of Mathematics\\ Hangzhou Normal University\\ Hang -zhou, China}
\email{<huanyinchen@aliyun.com>}
\address{Farzanegan Campus, Semnan University, Semnan, Iran}
\email{<m.sheibani@semnan.ac.ir>}

 \thanks{$^*$Corresponding author}

\subjclass[2010]{15A09, 15A10.} \keywords{Drazin inverse; additive formula; spectral idempotent; block matrix.}

\begin{abstract}
We present a new formula for the Drazin inverse of the sum of two complex matrices under weaker conditions with perturbations. By using this additive results, we establish new representations for the Drazin inverse of $2\times 2$ block complex matrix $\left(
\begin{array}{cc}
A&B\\
C&D
\end{array}
\right)$ under certain new conditions involving perturbations. These extend many known results, e.g., Yang and Liu (J. Comput. Applied Math., 235(2011), 1412--1417), Dopazo and Martinez-Serrano (Linear Algebra Appl., 432(2010), 1896--1904).
\end{abstract}

\maketitle

\section{Introduction}

Let ${\Bbb C}^{n\times n}$ be the set of all $n\times n$ matrices over the complex field ${\Bbb C}$, and let $A\in {\Bbb C}^{n\times n}$. The Drazin inverse of $A$ is the unique matrix $A^D\in {\Bbb C}^{n\times n}$ satisfying the following relations: $$AA^D=A^DA, A^D=A^DAA^D~\mbox{and} ~A^m=A^{m+1}A^D$$ for some $m\in {\Bbb N}$. As is well known, $A$ has Drazin inverse if and only if $rank(A^m)=rank(A^{m+1})$ for some $m\in {\Bbb N}$. The Drazin invertibility of a complex matrix is attractive. It has widespread applications in singular differential equations, Markov chains, and iterative methods (see~\cite{B}). Many authors have studied such problems from many different views, e.g., ~\cite{B,S,Y,Y2,Z1,ZM}.

Let $P,Q\in {\Bbb C}^{n\times n}$. It is attractive to give the explicit formula of $(P+Q)^D$. The formula was firstly obtained by Drazin in the case of $PQ=QP=0$.
In 2001, Hartwig et al. gave the formula when $PQ=0$ in ~\cite{H}. In 2011, the additive formula was given under the new condition
$PQ^2=0, PQP=0$ (see ~\cite[Theorem 2.1]{Y}). This is an elementary result and it deduce many wider conditions under which the Drazin inverse is expressed, e.g., ~\cite{S,Y2}. We denote by $A^{\pi}$ the eigenprojection of $A$ corresponding to the eigenvalue $0$ that is given by $A^{\pi}=I-AA^D$. The aim of this paper is to generalize elementary result. In Section 2, we extend ~\cite[Theorem 2.1]{Y} and present a representation of $(P+Q)^D$ in the case of $PQ^2=0$ and $PQP(PQ)^{\pi}=0$.

Let $M=\left(
  \begin{array}{cc}
    A&B\\
    C&D
  \end{array}
\right)$, where $A\in {\Bbb C}^{m\times m}, B\in {\Bbb C}^{m\times n}, C\in {\Bbb C}^{n\times m},$ $D\in {\Bbb C}^{n\times n}$. It is of interesting to find the Drazin inverse of the block complex matrix $M$. This problem is quite complicated and was expensively studied by many authors, see for example~\cite{D1,D,Y,Y2}. In Section 3, we apply these computational formulas to give the Drazin inverse of a block complex matrix $M$. If $BCB=0, DCB=0, BCA(BC)^{\pi}=0$ and $DCA(BC)^{\pi}=0$, we establish the representation of $M^D$, which also extend the result of Yang and Liu (see~\cite[Theorem 3.1]{Y}).

\section{additive results}

In this section, we will give a new additive formula for the Drazin inverse under a weaker condition, for which the following has been developed.
We begin with

\begin{lem} Let $P,Q\in {\Bbb C}^{n\times n}$. If $PQ=0$, then $$(P+Q)^D=\sum\limits_{i=0}^{\infty}(Q^D)^{i+1}P^iP^{\pi}+\sum\limits_{i=0}^{\infty}Q^iQ^{\pi}(P^D)^{i+1}.$$\end{lem}
\begin{proof} See~\cite[Lemma 1.2]{Y2}.\end{proof}

\begin{lem} Let $M=\left(
\begin{array}{cc}
E&I\\
F&0
\end{array}
\right), E,F\in {\Bbb C}^{n\times n}$. If $FEF^{\pi}=0$, then $$\begin{array}{lll}
M^{D}&=&\left(
  \begin{array}{cc}
    \Gamma&\Delta\\
    \Lambda&\Xi\\
     \end{array}
\right)~~~~~~~~~~~~~~~~~~~~~~~~~~~~~~~~~~~ (*)
\end{array}$$
where $$\begin{array}{lll}
\Gamma&=&\sum\limits_{i=0}^{\infty}(E^D)^{2i+1}F^iF^{\pi}+(E^{\pi}-\sum\limits_{j=1}^{\infty}(E^D)^{2j}F^j)F^{\pi}EF^D\\
&+&\sum\limits_{i=1}^{\infty}[(E^{\pi}-\sum\limits_{j=1}^{\infty}(E^D)^{2j}F^j)F^iF^{\pi}E+E^{\pi}E^{2i}F^{\pi}EFF^D\\
&+&\sum\limits_{k=0}^{i-2}E^{\pi}E^{2k+2}F^{i-k-1}F^{\pi}E]F^DC_i+\sum\limits_{i=1}^{\infty}[E^{\pi}E^{2i-1}F^{\pi}EFF^D\\
&-&\sum\limits_{j=0}^{\infty}(E^D)^{2j+1}F^{i+j}F^{\pi}E+\sum\limits_{k=0}^{i-2}E^{\pi}E^{2k+1}F^{i-k-1}F^{\pi}E]A_i,\\
\Delta&=&-F^D-E^DF^{\pi}EF^D-\sum\limits_{i=1}^{\infty}(E^D)^{2i+1}F^iF^{\pi}EF^D+\sum\limits_{i=0}^{\infty}(E^D)^{2i}F^iF^{\pi}\\
&-&[E^{\pi}-\sum\limits_{j=1}^{\infty}(E^D)^{2j}F^j]F^{\pi}EF^DEF^D+\sum\limits_{i=1}^{\infty}[(E^{\pi}-\sum\limits_{j=1}^{\infty}(E^D)^{2j}F^j)F^iF^{\pi}E\\
&+&E^{\pi}E^{2i}F^{\pi}EFF^D+\sum\limits_{k=0}^{i-2}E^{\pi}E^{2k+2}F^{i-k-1}F^{\pi}E]F^DD_i\\
&+&\sum\limits_{i=1}^{\infty}[E^{\pi}E^{2i-1}F^{\pi}EFF^D-\sum\limits_{j=0}^{\infty}(E^D)^{2j+1}F^{i+j}F^{\pi}E\\
&+&\sum\limits_{k=0}^{i-2}E^{\pi}E^{2k+1}F^{i-k-1}F^{\pi}E]B_i,\\
\Lambda&=&FF^D+\sum\limits_{i=1}^{\infty}F^iF^{\pi}EA_i,\\
\Xi&=&-FF^DEF^D+\sum\limits_{i=1}^{\infty}F^iF^{\pi}EB_;\\
\end{array}$$
$$\begin{array}{lll}
A_1&=&-F^DE,\\
B_1&=&(F^D)^2+F^DE^2,\\
C_1&=&F^D+FF^DE^2,\\
D_1&=&-F^DE-FF^DEF^D-FF^DE^3;
\end{array}$$
$$\begin{array}{lll}
A_{i+1}&=&F^DA_i-F^DEC_i,\\
B_{i+1}&=&F^DB_i-F^DED_i,\\
C_{i+1}&=&-FF^DEA_i+(F^D+FF^DE^2)C_i,\\
D_{i+1}&=&-FF^DEB_i+(F^D+FF^DE^2)D_i.
\end{array}$$
\end{lem}
\begin{proof} Construct $A_i,B_i,C_i,D_i (i\in {\Bbb N})$ as in the preceding, we easily check that $$\begin{array}{rll}
\left(
\begin{array}{cc}
0&F^D\\
FF^D&-FF^DEF^D
\end{array}
\right)^{2i+1}&=&\left(
\begin{array}{cc}
A_i&B_i\\
C_i&D_i
\end{array}
\right),
\end{array}$$
$$\begin{array}{ll}
&\left(
\begin{array}{cc}
0&F^D\\
FF^D&-FF^DEF^D
\end{array}
\right)^{2i+2}\\
=&\left(
\begin{array}{cc}
0&F^D\\
FF^D&-FF^DEF^D
\end{array}
\right)\left(
\begin{array}{cc}
A_i&B_i\\
C_i&D_i
\end{array}
\right)\\
=&\left(
\begin{array}{cc}
F^DC_i&F^DD_i\\
FF^DA_i-FF^DEF^DC_i&FF^DB_i-FF^DEF^DD_i
\end{array}
\right).
\end{array}$$
As the g-Drazin and Drazin inverse for a complex matrix coincide with each other, we obtain the result by~\cite[Theorem 2.3]{Z3}.\end{proof}

We now apply the preceding lemmas to obtain the following useful result.

\begin{thm} Let $P,Q\in {\Bbb C}^{n\times n}$. If $PQ^2=0$ and $PQP(PQ)^{\pi}=0$, then $$(P+Q)^D=(I,Q)\big[\left(
\begin{array}{cc}
P&PQ\\
I&Q
\end{array}
\right)^D\big]^2\left(
\begin{array}{c}
P\\
I
\end{array}
\right),$$
where $$\left(
\begin{array}{cc}
P&PQ\\
I&Q
\end{array}
\right)^D=\sum\limits_{i=0}^{\infty}(L^D)^{i+1}K^iK^{\pi}+\sum\limits_{i=0}^{\infty}L^iL^{\pi}(K^D)^{i+1},$$
$$\begin{array}{c}
L=\left(
\begin{array}{cc}
0&0\\
0&Q
\end{array}
\right), L^D=\left(
\begin{array}{cc}
0&0\\
0&Q^D
\end{array}
\right),\\
K=\left(
\begin{array}{cc}
P&PQ\\
I&0
\end{array}
\right), K^D=\left(
\begin{array}{cc}
P&I\\
I&0
\end{array}
\right)(H^D)^2\left(
\begin{array}{cc}
I&0\\
0&PQ
\end{array}
\right),
\end{array}$$
$$H^D=\left(
  \begin{array}{cc}
    \Gamma&\Delta\\
    \Lambda&\Xi\\
     \end{array}
\right), E=P, F=PQ,$$ where $\Gamma, \Delta,\Lambda$ and $\Xi$ constructed as in $(*)$.
\end{thm}
\begin{proof} Clearly, we have $P+Q=(I,Q)\left(
\begin{array}{c}
P\\
I
\end{array}
\right)$. By using the Cline's formula (see~\cite[Lemma 1.5]{X}), we have
$$(P+Q)^D=(I,Q)\big[\left(
\begin{array}{cc}
P&PQ\\
I&Q
\end{array}
\right)^D\big]^2\left(
\begin{array}{c}
P\\
I
\end{array}
\right).$$
Moreover, we have $$\left(
\begin{array}{cc}
P&PQ\\
I&Q
\end{array}
\right)=K+L,$$ where $$K=\left(
\begin{array}{cc}
P&PQ\\
I&0
\end{array}
\right), L=\left(
\begin{array}{cc}
0&0\\
0&Q
\end{array}
\right).$$
Let $$H=\left(
\begin{array}{cc}
P&I\\
PQ&0
\end{array}
\right), E=P, F=PQ.$$ Then $FEF^{\pi}=PQP(PQ)^{\pi}=0$. In view of
Lemma 2.2, $$H^D=\left(
  \begin{array}{cc}
    \Gamma&\Delta\\
    \Lambda&\Xi\\
     \end{array}
\right),$$ where $\Gamma, \Delta,\Lambda$ and $\Xi$ are constructed as in $(*)$.

Obviously, we check that
$$\begin{array}{lll}
K&=&\left(
\begin{array}{cc}
P&I\\
I&0
\end{array}
\right)\left(
\begin{array}{cc}
I&0\\
0&PQ
\end{array}
\right),\\
H&=&\left(
\begin{array}{cc}
I&0\\
0&PQ
\end{array}
\right)\left(
\begin{array}{cc}
P&I\\
I&0
\end{array}
\right).
\end{array}$$
By using Cline's formula, we have
$$K^D=\left(
\begin{array}{cc}
P&I\\
I&0
\end{array}
\right)(H^D)^2\left(
\begin{array}{cc}
I&0\\
0&PQ
\end{array}
\right).$$
Since $KL=0$, it follows by Lemma 2.1 that
$$\left(
\begin{array}{cc}
P&PQ\\
I&Q
\end{array}
\right)^D=\sum\limits_{i=0}^{\infty}(L^D)^{i+1}K^iK^{\pi}+\sum\limits_{i=0}^{\infty}L^iL^{\pi}(K^D)^{i+1},$$ where
$$\begin{array}{rll}
(K^D)^i&=&\left(
\begin{array}{cc}
P&I\\
I&0
\end{array}
\right)(H^D)^{i+1}\left(
\begin{array}{cc}
I&0\\
0&PQ
\end{array}
\right),\\
K^{\pi}&=&I_2-\left(
\begin{array}{cc}
P&I\\
I&0
\end{array}
\right)H^D\left(
\begin{array}{cc}
I&0\\
0&PQ
\end{array}
\right),\\
\end{array}$$
$$K^iK^{\pi}=\left(
\begin{array}{cc}
P&I\\
I&0
\end{array}
\right)H^{i-1}H^{\pi}\left(
\begin{array}{cc}
I&0\\
0&PQ
\end{array}
\right) (i\geq 1).$$ This completes the proof.\end{proof}

For future use, we now record the following.

\begin{cor} Let $P,Q\in {\Bbb C}^{n\times n}$. If $Q^2=0$ and $PQP(PQ)^{\pi}=0$, then $$\begin{array}{c}
(P+Q)^D=(I,Q)\big[K^D+L(K^D)^2\big]^2\left(
\begin{array}{c}
P\\
I
\end{array}
\right), E=P, F=PQ,\\
L=\left(
\begin{array}{cc}
0&0\\
0&Q
\end{array}
\right), K^D=\left(
\begin{array}{cc}
P&I\\
I&0
\end{array}
\right)\left(
  \begin{array}{cc}
    \Gamma&\Delta\\
    \Lambda&\Xi\\
     \end{array}
\right)^2\left(
\begin{array}{cc}
I&0\\
0&PQ
\end{array}
\right),
\end{array}$$ where $\Gamma, \Delta,\Lambda$ and $\Xi$ are constructed as in $(*)$.
\end{cor}
\begin{proof} Since $KL=0$ and $L^2=0$, it follows by Theorem ??? that
$$\begin{array}{lll}
\left(
\begin{array}{cc}
P&PQ\\
I&Q
\end{array}
\right)^D&=&(K+L)^D\\
&=&K^D+L(K^D)^2.
\end{array}$$ Therefore we obtain the result by Theorem 2.3.\end{proof}

As a direct consequence of Theorem 2.3, we derive

\begin{cor} (see~\cite[Theorem 2.1]{Y}) Let $P,Q\in {\Bbb C}^{n\times n}$. If $PQ^2=0$ and $PQP=0$, then $$\begin{array}{lll}
(P+Q)^D&=&\sum\limits_{i=0}^{\infty}(Q^D)^{i+1}P^iP^{\pi}+\sum\limits_{i=0}^{\infty}Q^iQ^{\pi}(P^D)^{i+1}\\
&+&\sum\limits_{i=0}^{\infty}Q^iQ^{\pi}(P^D)^{i+2}Q+\sum\limits_{i=0}^{\infty}(Q^D)^{i+3}P^{i+1}P^{\pi}\\
&-&Q^DP^DQ-(Q^D)^2PP^DQ.
\end{array}$$\end{cor}

\section{block complex matrices}

The purpose of this section is to use the preceding additive results to give some representations for the Drazin inverse of block matrix $M$.
We come now to extend ~\cite[Theorem 3.1]{Y} as follows.

\begin{thm} Let $M=\left(
  \begin{array}{cc}
    A & B \\
    C & D
  \end{array}
\right)\in {\Bbb C}^{(m+n)\times (m+n)}$. If $BCB=0, DCB=0, BCA(BC)^{\pi}=0$ and $DCA(BC)^{\pi}=0$, then $$\begin{array}{c}
M^D=\left(
\begin{array}{cccc}
I&0&0&0\\
0&I&C&0
\end{array}
\right)
\big[K^D+L(K^D)^2\big]^2
\left(
\begin{array}{cc}
A&B\\
0&D\\
I&0\\
0&I
\end{array}
\right),\\
E=\left(
\begin{array}{cc}
A&B\\
0&D
\end{array}
\right), F=\left(
\begin{array}{cc}
BC&0\\
DC&0
\end{array}
\right),\\
L=\left(
\begin{array}{cccc}
0&0&0&0\\
0&0&0&0\\
0&0&0&0\\
0&0&C&0
\end{array}
\right),\\
K^D=\left(
\begin{array}{cccc}
A&B&I&0\\
0&D&0&I\\
I&0&0&0\\
0&I&0&0
\end{array}
\right)\left(
  \begin{array}{cc}
    \Gamma&\Delta\\
    \Lambda&\Xi\\
     \end{array}
\right)^2\left(
\begin{array}{cccc}
I&0&0&0\\
0&I&0&0\\
0&0&BC&0\\
0&0&DC&0
\end{array}
\right),
\end{array}$$ where $\Gamma, \Delta,\Lambda$ and $\Xi$ are constructed as in $(*)$.
\end{thm}
\begin{proof} Clearly, we have $M=P+Q$, where
$$P=\left(
\begin{array}{cc}
A&B\\
0&D
\end{array}
\right), Q=\left(
\begin{array}{cc}
0&0\\
C&0
\end{array}
\right).$$ Then
$$Q^2=0, PQP(PQ)^{\pi}=0.$$
Therefore we complete the proof by Corollary 2.4.\end{proof}

As a consequence, we can extend~\cite[Theorem 2.2]{D} as follows:

\begin{cor} Let $M=\left(
  \begin{array}{cc}
    A & B \\
    C & D
  \end{array}
\right)\in {\Bbb C}^{(m+n)\times (m+n)}$. If $BCB=0, BDC=0, BD^2=0$ and $BCA(BC)^{\pi}=0$, then
$$\begin{array}{rll}
M^D&=&Q^{\pi}\sum\limits_{i=0}^{\infty}Q^i(P^D)^{i+1}+\sum\limits_{i=0}^{\infty}(Q^D)^{i+1}P^iP^{\pi}+Q^{\pi}\sum\limits_{i=0}^{\infty}Q^i(P^D)^{i+2}Q\\
&+&\sum\limits_{i=0}^{\infty}(Q^D)^{i+3}P^{i+1}P^{\pi}Q-Q^DP^DQ-(Q^D)^2PP^DQ,\\
\end{array}$$
$$\begin{array}{c}
P=\left(
  \begin{array}{cc}
    A & B \\
    C & 0
  \end{array}
\right), Q=\left(
  \begin{array}{cc}
    0 & 0 \\
    0 & D
  \end{array}
\right),\\
P^D=\left(
\begin{array}{cccc}
I&0&0&0\\
0&I&C&0
\end{array}
\right)
\big[K^D+L(K^D)^2\big]^2
\left(
\begin{array}{cc}
A&B\\
0&0\\
I&0\\
0&I
\end{array}
\right),\\
E=\left(
\begin{array}{cc}
A&B\\
0&0
\end{array}
\right), F=\left(
\begin{array}{cc}
BC&0\\
0&0
\end{array}
\right),\\
L=\left(
\begin{array}{cccc}
0&0&0&0\\
0&0&0&0\\
0&0&0&0\\
0&0&C&0
\end{array}
\right),\\
K^D=\left(
\begin{array}{cccc}
A&B&I&0\\
0&0&0&I\\
I&0&0&0\\
0&I&0&0
\end{array}
\right)\left(
  \begin{array}{cc}
    \Gamma&\Delta\\
    \Lambda&\Xi\\
     \end{array}
\right)^2\left(
\begin{array}{cccc}
I&0&0&0\\
0&I&0&0\\
0&0&BC&0\\
0&0&0&0
\end{array}
\right),
\end{array}$$ where $\Gamma, \Delta,\Lambda$ and $\Xi$ are constructed as in $(*)$.\end{cor}\begin{proof} Let $$P=\left(
\begin{array}{cc}
A&B\\
C&0
\end{array}
\right), Q=
\left(
\begin{array}{cc}
0&0\\
0&D
\end{array}
\right).$$ Since $BDC=0$ and $BD^2=0$, we verify that
$PQ^2=0$ and $PQP=0$. In view of Corollary 2.5,
$$\begin{array}{rll}
M^D&=&Q^{\pi}\sum\limits_{i=0}^{\infty}Q^i(P^D)^{i+1}+\sum\limits_{i=0}^{\infty}(Q^D)^{i+1}P^iP^{\pi}\\
&+&Q^{\pi}\sum\limits_{i=0}^{\infty}Q^i(P^D)^{i+2}Q+\sum\limits_{i=0}^{\infty}(Q^D)^{i+3}P^{i+1}P^{\pi}Q\\
&-&Q^DP^DQ-(Q^D)^2PP^DQ,\\
\end{array}$$ This complete the proof by applying Theorem 3.1 to $P$.\end{proof}

The following example illustrates that Corollary 3.4 is a nontrivial generalization of~\cite[Theorem 2.2]{D}.

\begin{exam} Let $M=\left(
  \begin{array}{cc}
    A & B \\
    C & D
  \end{array}
\right)\in M_4({\Bbb C})$, where $$A=\left(
  \begin{array}{cc}
    0 & 1 \\
    0 & 0
  \end{array}
\right), B=\left(
  \begin{array}{cc}
    0 & 2 \\
    0 & 0
  \end{array}
\right), C=\left(
  \begin{array}{cc}
    0 & 1 \\
    0 & -1
  \end{array}
\right), D=\left(
  \begin{array}{cc}
    -1 & 1 \\
    0 & 0
  \end{array}
\right).$$
Then $BCB=0, BDC=0, BD^2=0$ and $BCA(BC)^{\pi}=0$ while $BC\neq 0$. Clearly, we compute that
$$A^D=0, D^D=\left(
  \begin{array}{cc}
    -1 & 1 \\
    0 & 0
  \end{array}
\right).$$ Construct $E,F,P$ and $Q$ as in Theorem 3.1. Obviously, we have $Q^2=0$. Hence, by Corollary 3.2, we obtain the exact value $M^D$:
$$\begin{array}{lll}
M^D&=&P^D+Q(P^D)^2\\
&=&\left(
  \begin{array}{cccc}
    0&0&0&0\\
    0&0&0&0\\
    0&2&-1&1\\
    0&0&0&0\\
  \end{array}
\right).
\end{array}$$ \end{exam}

In a similar way as it was done in Theorem 3.1, using the another splitting, we have

\begin{thm} Let $M=\left(
  \begin{array}{cc}
    A & B \\
    C & D
  \end{array}
\right)\in {\Bbb C}^{(m+n)\times (m+n)}$. If $ABC=0, CBC=0, ABD(CB)^{\pi}=0$ and $CBD(CB)^{\pi}=0$, then $$\begin{array}{c}
M^D=\left(
\begin{array}{cccc}
I&0&0&B\\
0&I&0&0
\end{array}
\right)\big[K^D+L(K^D)^2\big]^2\left(
\begin{array}{cc}
A&0\\
C&D\\
I&0\\
0&I
\end{array}
\right),\\
E=\left(
\begin{array}{cc}
A&0\\
C&D
\end{array}
\right), F=\left(
\begin{array}{cc}
0&AB\\
0&CB
\end{array}
\right),\\
L=\left(
\begin{array}{cccc}
0&0&0&0\\
0&0&0&0\\
0&0&0&B\\
0&0&0&0
\end{array}
\right),\\
K^D=\left(
\begin{array}{cccc}
A&0&I&0\\
C&D&0&I\\
I&0&0&0\\
0&I&0&0
\end{array}
\right)\left(
  \begin{array}{cc}
    \Gamma&\Delta\\
    \Lambda&\Xi\\
     \end{array}
\right)^2\left(
\begin{array}{cccc}
I&0&0&0\\
0&I&0&0\\
0&0&0&AB\\
0&0&0&CB
\end{array}
\right),
\end{array}$$ where $\Gamma, \Delta,\Lambda$ and $\Xi$ are constructed as in $(*)$.
\end{thm}
\begin{proof} Write $M=P+Q$,  where
$$P=\left(
\begin{array}{cc}
A&0\\
C&D
\end{array}
\right), Q=\left(
\begin{array}{cc}
0&B\\
0&0
\end{array}
\right).$$ Then one checks that $$Q^2=0, PQP(PQ)^{\pi}=0.$$
Therefore we complete the proof by using Corollary 2.4.\end{proof}

As a consequence of the preceding theorem, we now derive

\begin{cor}Let $M=\left(
  \begin{array}{cc}
    A & B \\
    C & D
  \end{array}
\right)\in {\Bbb C}^{(m+n)\times (m+n)}$. If $CAB=0, CBC=0, A^2B=0$ and $CBD(CB)^{\pi}=0$, then $$\begin{array}{rll}
M^D&=&Q^{\pi}\sum\limits_{i=0}^{\infty}Q^i(P^D)^{i+1}+\sum\limits_{i=0}^{\infty}(Q^D)^{i+1}P^iP^{\pi}+P\sum\limits_{i=0}^{\infty}(Q^D)^{i+2}P^iP^{\pi}\\
&+&P\sum\limits_{i=0}^{\infty}Q^{\pi}Q^{i+1}(P^D)^{i+3}-PQ^DP^D-PQQ^D(P^D)^2,\\
\end{array}$$
$$\begin{array}{c}
P=\left(
  \begin{array}{cc}
    A & 0 \\
    0 & 0
  \end{array}
\right), Q=\left(
  \begin{array}{cc}
    0 & B \\
    C & D
  \end{array}
\right),
\end{array}$$
$$\begin{array}{c}
Q^D=\left(
\begin{array}{cccc}
I&0&0&B\\
0&I&0&0
\end{array}
\right)\big[K^D+L(K^D)^2\big]^2\left(
\begin{array}{cc}
0&0\\
C&D\\
I&0\\
0&I
\end{array}
\right),\\
E=\left(
\begin{array}{cc}
0&0\\
C&D
\end{array}
\right), F=\left(
\begin{array}{cc}
0&0\\
0&CB
\end{array}
\right),\\
L=\left(
\begin{array}{cccc}
0&0&0&0\\
0&0&0&0\\
0&0&0&B\\
0&0&0&0
\end{array}
\right),\\
K^D=\left(
\begin{array}{cccc}
0&0&I&0\\
C&D&0&I\\
I&0&0&0\\
0&I&0&0
\end{array}
\right)\left(
  \begin{array}{cc}
    \Gamma&\Delta\\
    \Lambda&\Xi\\
     \end{array}
\right)^2\left(
\begin{array}{cccc}
I&0&0&0\\
0&I&0&0\\
0&0&0&0\\
0&0&0&CB
\end{array}
\right),
\end{array}$$ where $\Gamma, \Delta,\Lambda$ and $\Xi$ are constructed as in $(*)$.
\end{cor}
\begin{proof} Write $M=P+Q$,  where
$$P=\left(
\begin{array}{cc}
A&0\\
0&0
\end{array}
\right), Q=\left(
  \begin{array}{cc}
    0 & B \\
    C & D
  \end{array}
\right).$$ Then we have $$P^2Q=0, QPQ=0.$$
In view of \cite[Theorem 2.2]{Y}, we have
$$\begin{array}{lll}
M^D&=&Q^{\pi}\sum\limits_{i=0}^{\infty}Q^i(P^D)^{i+1}+\sum\limits_{i=0}^{\infty}(Q^D)^{i+1}P^iP^{\pi}+P\sum\limits_{i=0}^{\infty}(Q^D)^{i+2}P^iP^{\pi}\\
&+&P\sum\limits_{i=0}^{\infty}Q^{\pi}Q^{i+1}(P^D)^{i+3}-PQ^DP^D-PQQ^D(P^D)^2.
\end{array}$$ Applying Theorem 2.3 to $Q$, we complete the proof.\end{proof}

\end{document}